\documentclass[a4paper]{amsart}
\usepackage{graphics}
\usepackage{amsmath}
\usepackage{latexsym}
\usepackage{amssymb}
\usepackage[all]{xy}
\xyoption{matrix}
\xyoption{arrow}

\begin{document}

\newtheorem{lem}{Lemma}[section]
\newtheorem{prop}[lem]{Proposition}
\newtheorem{cor}[lem]{Corollary}
\newtheorem{thm}[lem]{Theorem}
\newtheorem*{thmA}{Theorem}
\newtheorem*{thmB}{Theorem B}
\newtheorem{rem}[lem]{Remark}
\newtheorem{defin}[lem]{Definition}

\newcommand{\Z}{\operatorname{\mathbb Z}\nolimits}
\newcommand{\Q}{\operatorname{\mathbb Q}\nolimits}
\newcommand{\Ext}{\operatorname{Ext}\nolimits}
\newcommand{\Hom}{\operatorname{Hom}\nolimits}
\newcommand{\add}{\operatorname{add}\nolimits}
\renewcommand{\mod}{\operatorname{mod}\nolimits}
\renewcommand{\dim}{\operatorname{dim}\nolimits}
\newcommand{\End}{\operatorname{End}\nolimits}
\newcommand{\op}{\operatorname{op}\nolimits}
\newcommand{\C}{\operatorname{\mathcal C}\nolimits}
\newcommand{\D}{\operatorname{\mathcal D}\nolimits}
\newcommand{\M}{\operatorname{\mathcal M}\nolimits}
\renewcommand{\P}{\operatorname{\mathcal P}\nolimits}
\renewcommand{\L}{\Lambda}
\newcommand{\G}{\Gamma}

\title[Tilted versus cluster-tilted]
{From tilted to cluster-tilted algebras of Dynkin type}

\author[Buan]{Aslak Bakke Buan}
\address{Institutt for matematiske fag\\
Norges teknisk-naturvitenskapelige universitet\\
N-7491 Trondheim\\
Norway}
\email{aslakb@math.ntnu.no}

\author[Reiten]{Idun Reiten}
\address{Institutt for matematiske fag\\
Norges teknisk-naturvitenskapelige universitet\\
N-7491 Trondheim\\
Norway}
\email{idunr@math.ntnu.no}

\begin{abstract}
We show how a cluster-tilted algebra of finite representation type is related
to the corresponding tilted algebra, in the case
of algebras defined over an algebraically closed field. 
\end{abstract}

\maketitle

\section{Introduction}

The cluster categories $\C_H$ of finite dimensional hereditary algebras $H$ were
introduced and studied in \cite{bmrrt} in order to give a categorical model 
for some of the central ingredients of the theory of cluster algebras introduced by Fomin and Zelevinsky
\cite{fz}.
In particular tilting theory for cluster categories was investigated.
Associated with a tilting module $T$ over a hereditary algebras is, in addition to the 
tilted algebra $\L = \End_{H}(T)^{\op}$, also the cluster-tilted algebra  $\G = \End_{\C_H}(T)^{\op}$,
investigated in \cite{bmr1, bmr2, bmr3}.

It is of interest to compare the algebras $\L$ and $\G$ associated with the same tilting module. In general
we do not know how to construct $\G$ from $\L$. Here we deal with 
the case where $H$ is of Dynkin type over an algebraically closed
field $K$. In this case we show how the cluster-tilted algebra $\G$
is determined by the quiver and relations of the tilted algebra $\L$.

\section{Background}

In this section we recall some notions and basic results.
Let $K$ be an algebraically closed field, and $Q$ a finite quiver with no oriented cycles.
Then the path algebra $H = KQ$ is a hereditary finite-dimensional algebra. Let $\mod H$ be
the category of (left) finitely generated modules over $H$
and let $\underline{r}$ denote the Jacobson radical of an algebra.

\subsection{Tilting modules and tilted algebras}

A module $T$ in $\mod H$ is called a {\em tilting module} if $\Ext^1_H(T,T) = 0$ and if $T$ is maximal
with respect to this property, i.e. if $\Ext^1_H(T \amalg X,T \amalg X) = 0$, then
$X$ is a direct summand in a direct sum of copies of $T$.
The endomorphism-algebras of tilting modules are called {\em tilted algebras}.
See \cite{h} for basic properties of tilted algebras and further references.

\subsection{Mutation of quivers}

In connection with their definition of cluster algebras,
Fomin and Zelevinsky \cite{fz} defined a {\em mutation}-operation on skew-symmetrizable 
integer matrices.
Mutation of quivers can be seen as a special case of this, and can be
defined for any finite quiver with no loops and no oriented cycles of length two. 

A quiver is called {\em double path avoiding} if it does not contain multiple arrows,
and if any quiver obtained by repeated mutating and/or factoring (i.e. removing one or more vertices)
also contains no multiple arrows. 
See Section 1 of \cite{bmr3} for more details on this.
A fact that we will freely use here is that
the quiver of a cluster-tilted algebra of finite representation type is double path avoiding.

\subsection{Cluster categories and cluster-tilted algebras}

Cluster categories were introduced in \cite{bmrrt} in order to give a categorical
approach to cluster algebras as defined by Fomin and Zelevinsky \cite{fz}. In 
\cite{ccs1} a category was introduced for Dynkin quivers of type $A_n$, which was shown
to be equivalent to the cluster category. 
The idea was to model the combinatorics of clusters on the set of tilting 
objects in a cluster category.
The first link from cluster algebras to tilting theory was given by
Marsh, Reineke and Zelevinsky \cite{mrz}.
Given a hereditary algebra $H$, the cluster category is defined as follows.
Let $\D = D^b(\mod H)$ be the bounded derived category of $H$. This category has a
triangulated structure. Let $[1]$ denote the suspension functor. It also has AR-triangles, by
\cite{h}. Let $\tau$ denote the corresponding auto-equivalence on $\D$ with quasi-inverse $\tau^{-1}$. 
Consider the auto-equivalence
$F = \tau^{-1}[1]$. The cluster category $\C$ is the orbit category $\C/F$. That is,
the objects of $\C$ are the objects of $\D$, while $\Hom_{\C}(X,Y) = \amalg_i \Hom_{\D}(X, F^i Y)$.
The indecomposable objects of $\D$ are stalk-complexes, and we identify $\mod H$ with the full
subcategory of $\D$ where the objects are stalk-complexes in degree 0.
It was shown by Keller \cite{k} that the cluster category is also triangulated and that
the canonical functor $\D \to \C$ is a triangle functor.
It was shown in \cite{bmrrt} that $\C$ is a Krull-Schmidt category with 
AR-triangles, and that the canonical functor $\D \to \C$ preserves AR-triangles.
    
An object $T$ of $\C$ is called a (cluster-)tilting object if $\Ext^1_{\C}(T,T) = 0$ and $T$ is maximal
with respect to this property, i.e. if $\Ext^1_{\C}(T \amalg X,T \amalg X) = 0$, then
$X$ is a direct summand in a direct sum of copies of $T$.
The endomorphism-algebra $\End_{\C}(T)^{\op}$ of a tilting object $T$ is called a {\em cluster-tilted algebra}.
A cluster-tilted algebra coming from $H$ is of finite representation
type if and only if $H$ is of finite representation type \cite{bmr1}.

We also need the following results from \cite{bmrrt}, \cite{bmr1} and \cite{bmr2} 
about tilting objects and cluster-tilted algebras.

\begin{thm}\label{gb}
Let $\C = \C_H$ for a hereditary algebra $H$, and
let $T$ be a tilting object in $\C$.
\begin{itemize}
\item[(a)]{Let $Q$ be the quiver of a cluster-tilted algebra $\G = \End_{\C}(T)^{\op}$.
Let $Q'$ be the quiver obtained by mutating $Q$ at the vertex $k$. Furthermore, let $T_k$ be
the indecomposable direct summand of $T$ corresponding to the vertex $k$ and assume $T = \bar{T} \amalg T_k$. 
Then there is a tilting object $T' = \bar{T} \amalg T_k^{\ast}$, such that $Q'$ is the quiver of
$\End_{\C}(T')^{\op}$.
}
\item[(b)]
{Assume $H$ is of finite type, and let $T_a$ and $T_b$
be indecomposable direct summands in $T$. Then $\dim \Hom_{\C}(T_a,T_b) \leq 1$.
If there is a non-zero map $T_a \to T_b$ in $\C$, which lifts to a map $T_a \to FT_b$ in $\D$ with $T_a, T_b$
in $\mod H$,
then any non-zero map $T_b \to T_a$ in $\C$ lifts to a map $T_b \to T_a$ in $\D$, with $T_a, T_b$
in $\mod H$.
}
\end{itemize}
\end{thm}



\subsection{Cluster-tilted algebras of finite representation type}

The proof of the main result of this paper strongly depends on results in \cite{bmr3}.
In that paper, the main aim was to show that the cluster-tilted algebras of finite type are
determined by their quivers and to describe the relations.
Some of the results and notions needed to prove this
are also useful here, so we recall them.

An oriented cycle in a quiver is called {\em full} if there are no repeated vertices
and if the subquiver generated by the cycle
contains no further arrows.
If there is an arrow $i \to j$ in
a quiver $Q$, then a path from $j$ to $i$ is called {\em a shortest path} if the 
induced subquiver is a full cycle.
Let $\G = KQ/I$ be a cluster-tilted algebra.
The elements in $I$ are called {\em relations} if they are linear combinations 
$k_1 \omega_1 + \cdots + k_m \omega_m$ of paths $\omega_i$ in $Q$, all starting in the same vertex and
ending in the same vertex, and with each $k_i$ non-zero in $K$. 
If $m=1$, we call the relation a 
{\em zero-relation}. If $m=2$, we call it a {\em commutativity-relation} (and say that the paths $\omega_1$ and
$\omega_2$ commute).
A relation $\rho$ is called {\em minimal} if 
whenever $\rho = \sum_i \beta_i \circ \rho_i \circ \gamma_i$, where $\rho_i$ is a relation
for every $i$, then there is an index $j$ such that both $\beta_j$ and $\gamma_j$ are scalars.

For ease of notation we use the same symbol $\omega$ to denote a path, the corresponding element in the path algebra,
and the corresponding element in $KQ/I$.

The following was the main result of \cite{bmr3}. A consequence is that a
cluster-tilted algebra of finite representation type is determined by its quiver, up to isomorphism. 

\begin{thm}\label{gc}
Let $\G = KQ/I$ be a cluster-tilted algebra of finite representation type, and let $i,j$ 
be vertices in $Q$.
\begin{itemize}
\item[(a)]{The ideal $I$ is generated by minimal zero-relations and minimal commutativity relations.}
\item[(b)]{Assume there is an arrow $i \to j$. Then there are at most two shortest paths
from $j$ to $i$.
\begin{itemize}
\item[(i)]{If there is exactly one, then this is a minimal zero-relation.}
\item[(ii)]{If there are two, $\omega$ and $\mu$, then  $\omega$ and $\mu$ are not zero in $\G$,
and there is a minimal relation $\omega + \lambda \mu$ for some $\lambda \neq 0$ in $K$.
}
\end{itemize}
}
\item[(c)]{Up to multiplication by non-zero elements of $K$, there are no other minimal zero-relations 
or commutativity relations than the ones coming from (b).}
\end{itemize}
\end{thm}

We also need some results which were used to prove the above theorem in \cite{bmr3}.

\begin{prop}\label{help}
Let $\G = KQ/I$ be a cluster-tilted algebra of finite representation type, and let $i,j$ 
be vertices in $Q$.
\begin{itemize}
\item[(a)]{A relation $\rho$ is minimal if and only if $\rho = \rho' +\rho''$, where
$\rho'$ is a minimal commutativity relation or a minimal zero-relation, and $\rho''$ is any relation. 
\item[(b)]{Given an arrow $i \to j$, and two distinct shortest paths $\omega$ and $\mu$ from $j$ to $i$,
$\omega$ and $\mu$
are both non-zero, and there is a minimal commutativity relation
$\omega + \lambda \mu$ for some $\lambda$ in $K$, where $\lambda$ can be assumed to be $-1$.}
\item[(c)]{For any minimal commutativity relation involving paths $\omega$ and $\mu$, the subquiver
generated by these paths has the form
$$
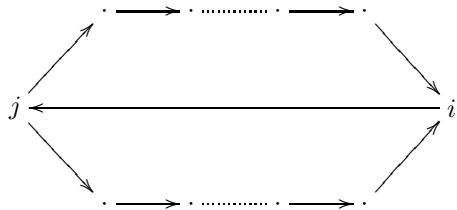
\xy
\xymatrix{
& \cdot \ar[r] & \cdot \ar@{.}[r] & \cdot \ar[r] & \cdot \ar[dr] & \\
j \ar[dr] \ar[ur] & & & & & i \ar[lllll] \\
& \cdot \ar[r] & \cdot \ar@{.}[r] & \cdot \ar[r] & \cdot \ar[ur] &
}
\endxy
$$
}
\item[(d)]{In a full cycle in $Q$ of length $n$, the composition of $n-2$ arrows is always non-zero.}
\item[(e)]{If the quiver $Q$ is an oriented cycle of length $n$, then $I = \underline{r}^{n-1}$.}
}
\end{itemize}
\end{prop}

\section{Main result}

Let $K$ denote an algebraically closed field,
and let $H$ be the path algebra over $K$ of some Dynkin quiver.
Let $T$ be a basic tilting module in $\mod H$, and let $\L = \End_H(T)^{\op}$ and $\G=  \End_{\C_H}(T)^{\op}$
be the corresponding tilted and cluster-tilted algebra, respectively.
Assume $\L = KQ'/I'$ and $\G = KQ/I$.
We will compare the quivers $Q$ and $Q'$.

Let $\add T$ denote the full subcategory of $\mod H$ with objects the direct summands of directs
sums of copies of $T$.
We let $(\add T)_{\C}$ denote the full subcategory of $\C = \C_H$ with objects the direct summands of directs
sums of copies of $T$.

The vertices $Q'_0$ of $Q'$ correspond to the indecomposable objects in $\add T$, while
the vertices $Q_0$ of $Q$ correspond to the indecomposable objects in $(\add T)_{\C}$.
It is clear that if we write $T$ as a direct sum of indecomposables 
$T = T_1 \amalg \dots \amalg T_n$
in $\mod H$, this is a also a decomposition into indecomposables in $\C$. We therefore
identify the sets $Q_0$ and $Q'_0$.

Let $T_a$ and $T_b$ be indecomposable objects of $\add T$. 
Using the definition of the cluster category $\C = \C_H$ and the fact that $\dim \Hom_{\C}(T_a,T_b) \leq 1$ by 
Theorem \ref{gb} (b), we have 
that any non-zero map $T_a \to T_b$ in $(\add T)_{\C}$ is either the image of a map $T_a \to T_b$ in $\D$, 
or the image of a map
$T_a \to FT_b$ in $\D$. 
Maps of the first kind are called {\em $m$-maps}, while maps of the second
kind are called {\em $f$-maps.}
Up to scalars, irreducible maps in $(\add T)_{\C}$ correspond to arrows in the quiver $Q$, while
irreducible maps in $\add T$ correspond to arrows in the quiver $Q'$. 
The arrows in $Q$ corresponding to irreducible $m$-maps, we call {\em $m$-arrows}, the other arrows
are called {\em $f$-arrows}.
For the set of arrows $Q_1$ in $Q$ we thus have a partition $Q_1 = Q_1^m \cup Q_1^f$, where
$Q_1^m$ are the $m$-arrows and $Q_1^f$ are the $f$-arrows.
A path in the quiver $Q$ is called an {\em $f$-path} if 
it contains at least one $f$-arrow. The other paths are called {\em $m$-paths}. 

\begin{lem}\label{preserved}
\begin{itemize}
\item[(a)]{An irreducible map $T_a \to T_b$ in $\add T$ is also irreducible 
as a map in $(\add T)_{\C}$.}
\item[(b)]{Every irreducible map $T_a \to T_b$ in $(\add T)_{\C}$ 
which is the image of a map in $\add T$ is also irreducible in $\add T$.}
\end{itemize}
\end{lem}

\begin{proof}
(a): Trivial, since an $m$-map in $\C_H$ can not factor through an $f$-map.

\noindent (b): Also trivial.
\end{proof}

This means that the set of arrows $Q'_1$ in $Q'$ corresponds to the set of $m$-arrows
$Q_1^m$ in $Q$. We want to show that the $f$-arrows are determined by the minimal relations
for $\L$. We first compare the minimal relations in $I' \subset KQ'$ with
the minimal relations in $I \subset KQ$. 

\begin{lem}
A minimal relation in $I' \subset KQ'$ is also a minimal relation in $I \subset KQ$.
\end{lem}

\begin{proof}
Given two indecomposable direct summands $T_a$ and $T_b$ in a tilting module $T$, we
have $\dim_k \Hom_{\C}(T_a, T_b) \leq 1$ by Theorem \ref{gb} (b). A consequence is that there cannot be
both non-zero $m$-maps and non-zero $f$-maps from $T_a$ to $T_b$.

Assume that $\rho'$ is a minimal relation in $I' \subset KQ'$. 
It is clear that $\rho'$ is also a relation in  
$I \subset KQ$. Since $I$ is generated by minimal commutativity-relations and minimal zero-relations,
by Theorem \ref{gc} (a), we can write $\rho' = \sum_{i=1}^m \beta_i \rho_i \gamma_i$ in $KQ$, with
$\rho_i$ minimal commutativity-relations for $i= 1, \dots, r$ and 
$\rho_i$ minimal zero-relations for $i= r+1, \dots, m$. Write $\rho_i$ as a sum of scalar multiples
of paths $\rho_i = \omega_i^{(0)} + \omega_i^{(1)}$, for
$i= 1, \dots, r$.
It follows from Theorem \ref{gb} (b)
that for $i= 1, \dots, r$ we have that $\omega_i^{(0)}$ is an $f$-path if and only
if $ \omega_i^{(1)}$ is an $f$-path. Hence, the terms containing $f$-paths in the sum 
$\sum_{i=1}^m \beta_i \rho_i \gamma_i$ cancel, and we are left with an expression of $\rho'$
as a linear combination of relations in $I' \subset KQ'$. By the minimality assumption on $\rho'$
in $I'$, this means that for some $i$ we have that 
$\beta_i$ and $\gamma_i$ are scalars. By Proposition \ref{help} (a), we have that $\rho'$ is also minimal in
$I \subset KQ$.
\end{proof}

Using Theorem \ref{gc} (c)
we obtain as a consequence of the above that for any minimal relation $\rho$ in $I' \subset KQ'$, there
is an $f$-arrow $\alpha_{\rho}$ in $Q$.
It remains to show that actually all $f$-arrows arise this way.
Recall that Theorem \ref{gc} (a) says that the the set of minimal zero-relations and minimal
commutativity relations in $I \subset KQ$ generate $I$.

\begin{thm}\label{main}
Let $\rho_1, \dots , \rho_d$ be the minimal relations in $I'$ which are either zero-relations or
commutativity relations.
Then the quiver $Q$ is given by $Q_0 = Q'_0$ and $Q_1 = Q_1' \cup \{\alpha_{\rho_i} \}$. 
\end{thm}

\begin{proof}
Let us fix some notation used in this proof.
Recall that if there is an arrow $i \to j$, then a path from $j$ to $i$ is called shortest if 
the induced oriented cycle is full. 
A {\em walk} in a quiver is a sequence of arrows 
$\alpha_1, \dots, \alpha_t$ such that $\alpha_j$ 
connects the vertices $i_j$ and $i_{j+1}$ for $j= 1, \dots, t-1$.
If all arrows involved are $m$-arrows it is called an $m$-walk.
A path in a quiver is hence a directed walk, and an $m$-path is a directed $m$-walk. 

It remains only to show that for every $f$-arrow, say $\alpha \colon j \to i$, there  
is a minimal relation involving paths from $i$ to $j$ in $I' \subset KQ'$.

We first show that there is an $m$-path from $i$ to $j$, which is a shortest path.
We show this by induction. 
Note that a walk of length two between $j$ and $i$ is necessarily an oriented path from $i$ to
$j$, since the quiver $Q$ is double path avoiding. 

We have the following observations.

\begin{lem}\label{four}
Every full oriented cycle of length $n \geq 4$ contains exactly one $f$-arrow.
\end{lem}

\begin{proof}
Since $Q'$ does not contain oriented cycles, it is clear that there is always at least one $f$-arrow
on every cycle. It is also clear that every path containing two (or more)
$f$-arrows must be zero, since such a path
corresponds to  a map $T_a \to F^2T_b$ in $\D$. On the other hand the composition of 
$n-2$ or less arrows on a cycle of length $n$ is non-zero by Proposition \ref{help} (d).
From this it follows that any full oriented cycle of length $\geq 5$ contains exactly one $f$-arrow.
It also follows that two consecutive $f$-arrows in a full oriented cycle of length 4 is impossible.
To exclude a full oriented cycle $c_1 \circ c_2 \circ c_3 \circ c_4$
where $c_1$ and $c_3$ are $f$-arrows and $c_2$ and $c_4$ are $m$-arrows, note that
in this case $c_1 \circ c_2$ is non-zero and $c_3 \circ c_4$ is non-zero, so we 
have a contradiction to Theorem \ref{gb} (b).
\end{proof}

It will actually follow from the remaining part of the proof, that also 
full oriented cycles of length 3 have this property.

\begin{lem}\label{commdia}
Assume $Q$ has a full subquiver of the form
$$
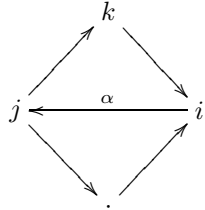
\xy
\xymatrix{
&  k \ar[dr] & \\
j  \ar[dr] \ar[ur] & & i \ar[ll]_{\alpha} \\
& \cdot \ar[ur] & 
}
\endxy
$$
\\
where $\alpha$ is an $f$-arrow. Then the other arrows in this full subquiver are all $m$-arrows.
\end{lem}

\begin{proof}
By Theorem \ref{gb} (b) and Proposition \ref{help} (c), it follows that either both
paths from $j$ to $i$ are $f$-paths, or both are $m$-paths.
Assume both are $f$-paths. Mutating at $k$ one obtains a full subquiver which is a 4-cycle. 
$$
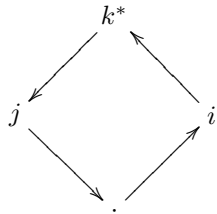
\xy
\xymatrix{
&  k^{\ast} \ar[dl] & \\
j  \ar[dr]  & & i \ar[ul] \\
& \cdot \ar[ur] & 
}
\endxy
$$
\\
Assume $\G = \End_{\C}(T)^{\op}$, for a tilting object $T$ in $\C_H$, and 
let $T_i, T_j$ be indecomposable direct summands in $T$ corresponding to the vertices $i$ and $j$, respectively.
Since $\dim \Hom_{\C}(T_i, T_j) \leq 1$, by Theorem \ref{gb} (b), 
and using Proposition \ref{help} (c), 
it is clear that
the arrow $\alpha$ and the composition $i \to k^{\ast} \to j$ represent the 
same map (up to scalars) in $(\add T)_{\C}$.
Thus the composition $i \to k^{\ast} \to j$ contains an $f$-arrow, and we have a contradiction to 
Lemma \ref{four}.
Thus all arrows except $\alpha$ are $m$-arrows.
\end{proof}

\begin{lem}\label{cycles}
Assume there is an arrow $x \to y$ in $Q$, such that there is a walk $\delta$ between $y$ to $x$.
Then there is a shortest path from $y$ to $x$ only passing through vertices on $\delta$.
\end{lem}

\begin{proof}
This follows from the fact that the quivers of cluster-tilted algebras of finite representation type
are double-path avoiding, and the fact \cite{bmr3} that any full subquiver which is a
non-oriented cycle, is not double path avoiding.
\end{proof}

The quiver $Q'$ is connected, since $H$ is connected. 
There is therefore always an $m$-walk between vertices
of $Q$. 
It is a direct consequence of Lemma \ref{cycles} that all $f$-arrows lie on oriented cycles.

Assume now that there is an $f$-arrow $\alpha \colon i \to j$ such that there
is no shortest $m$-path from $j$ to $i$. 
For each such $f$-arrow there is at least one $m$-walk between $i$  and $j$.
Let $m(\alpha)$ be the length of the $m$-walk between $i$ and $j$ of minimal length.
Amongst all $f$-arrows with the property that there is no $m$-path from $j$ to $i$, 
choose an $\alpha$ with smallest possible $m(\alpha)$.

By Lemma \ref{cycles} there is a shortest path from $j$ to $i$. By assumption this is not an
$m$-path.
Hence, by Lemma \ref{four}, it is of length two.
Assume the path is composed by the arrows $\beta \colon j \to k$ and $\gamma \colon k \to i$.
At least one of the arrows $\beta$ and $\gamma$ is an $f$-arrow. 
Consider also an $m$-walk $\delta$ between $i$ and $j$ of minimal length. Assume the walk 
is of length $t$, and consists of arrows $\delta_1, \dots, \delta_t$, where
$\delta_1$ is an arrow either starting or ending in $j$ and $\delta_t$ is an arrow either starting or ending in $i$.  
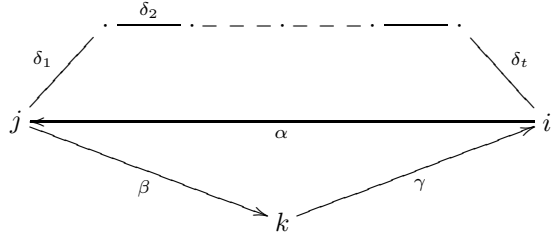
$$
\xy
\xymatrix{
& \cdot \ar@{-}[r]^{\delta_2} & \cdot \ar@{--}[r] & \cdot \ar@{--}[r] & \cdot \ar@{-}[r] & \cdot \ar@{-}[dr]^{\delta_t} & \\
j \ar[drrr]_{\beta} \ar@{-}[ur]^{\delta_1} & & & & & & i \ar[llllll]^{\alpha} \\
& & & k \ar[urrr]_{\gamma} & & &
}
\endxy
$$
\\
First assume
that $\gamma$ is an $f$-arrow, and that $\beta$ is an $m$-arrow. 
Consider the case where
$k$ is a vertex on the path $\delta$. Then by the assumption that $\delta$ is a walk
of minimal length, we must
have that $\beta = \delta_1 \colon j \to k$.
The induced $m$-walk $\delta'$ from $k$ to $i$ given by the arrows $\delta_2, \dots, \delta_t$
obviously has length smaller than the length of $\delta$, so since $\gamma$ is an $f$-arrow, and using the 
mimimality assumption, there exists a directed path of $m$-arrows from $i$ to $k$.
It is clear that $j$ is not a vertex on this path, so we have the following commutative
diagram.
$$
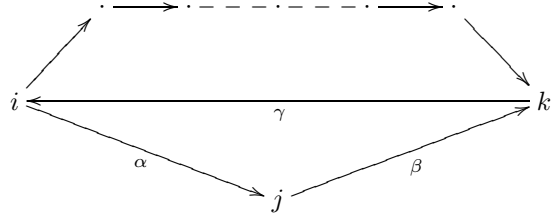
\xy
\xymatrix{
& \cdot \ar[r] & \cdot \ar@{--}[r] & \cdot \ar@{--}[r] & \cdot \ar[r] & \cdot \ar[dr] & \\
i \ar[drrr]_{\alpha} \ar[ur] & & & & & & k \ar[llllll]^{\gamma} \\
& & & j \ar[urrr]_{\beta} & & &
}
\endxy
$$
\\
By the minimality assumptions on the walk $\delta$ it is clear that
the two paths from $i$ to $k$ are such that the subquiver generated by these paths
is as in Proposition \ref{help} (c), and thus are non-zero and commute. 
By Theorem \ref{gb} (b) this is a contradiction, since one path is an $m$-path and 
the other is an $f$-path.

Now consider the case where $k$ is not on the walk $\delta$. Then, by Lemma \ref{cycles}, there is an alternative 
shortest path from $j$ to $i$ only passing through vertices on $\delta$. By assumption on $\alpha$
this path contains at least one $f$-arrow.
Thus, it is of length two. We get the following quiver
$$
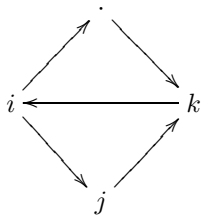
\xy
\xymatrix{
& \cdot \ar[dr] & \\
i \ar[ur] \ar[dr] & & k \ar[ll] \\
& j \ar[ur] & 
}
\endxy
$$
\\
with two shortest paths from $i$ to $k$ containing $f$-arrows,
and hence we have a contradiction to Lemma \ref{commdia}.
This finishes the case where $\gamma$ is an $f$-arrow, while $\beta$ is an $m$-arrow.

The case where $\beta$ is an $f$-arrow, while $\gamma$ is an $m$-arrow can be excluded by similar arguments.

Now assume that both $\beta$ and $\gamma$ are $f$-arrows.
Consider a shortest path $\omega$ from $j$ to $i$ induced by the minimal $m$-walk,
i.e. the path only passes through vertices on the chosen minimal $m$-walk. 
This path is necessarily of length two, by Lemma \ref{four}. There are two cases. 
We first consider the case that
$\omega$ is different from $\beta \circ \gamma$.
In this case, the subquiver generated
by these paths is as in Lemma \ref{commdia}, and 
hence we have a contradiction to this lemma.

The second case is when these paths are equal. Then consider the $f$-arrow $\beta \colon j \to k$.
The minimal $m$-walk between $i$ and $j$ induces a minimal $m$-walk $\delta''$ between
$j$ and $k$. By minimality, it is clear that $i$ is not a vertex on $\delta''$.
Thus the shortest path $\phi$ from $k$ to $j$ induced by $\delta''$ is disjoint from 
the path $\gamma \circ \alpha \colon k \to i \to j$. 

Since $\phi$ clearly
contains an $f$-arrow, it must be of length two. So, we have two different paths 
of length two from $k$ to $j$, and since $Q$ is double path avoiding it is clear that 
the subquiver generated by these paths is 
$$
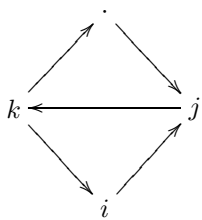
\xy
\xymatrix{
& \cdot \ar[dr] & \\
k \ar[ur] \ar[dr] & & j \ar[ll] \\
& i \ar[ur] & 
}
\endxy
$$
We thus have a contradiction by Lemma \ref{commdia}.

Let us now complete the proof of the theorem.
We have established that for every $f$-arrow $i \to j$, there is at least one shortest
path from $j$ to $i$ which is an $m$-path. By Theorem \ref{gc} there are either exactly one or exactly
two shortest paths from $j$ to $i$. From the same theorem it follows that in 
case there is exactly one, this is a minimal zero-relation.
Furthermore, in case there is a second shortest path from $j$ to $i$, this path must represent
a non-zero map, and hence is an $m$-path. 
Thus there is a minimal commutativity relation
involving these two paths, by Propostion \ref{help} (b).
Hence, in both cases there is a minimal relation in $I \subset KQ$ involving $m$-paths from $j$ to $i$.
It is clear that this is also a minimal relation in $I' \subset KQ'$.
\end{proof}

The following is obtained as a consequence of the last part of the proof of our main theorem,
using that there are no oriented cycles in the quiver of a tilted algebra.

\begin{cor}
Every full cycle in the quiver $Q$ of a cluster-tilted algebra 
of finite representation type contains exactly one $f$-map.
\end{cor}

As mentioned in the introduction, it is shown in \cite{bmr3} that
a cluster-tilted algebra of finite representation type is determined
by its quiver. 
This has the following consequence.

\begin{cor}
Let $H$ be the path algebra of some Dynkin quiver over an algebraically closed field $K$.
Given a tilted algebra $\L = \End_{H}(T)^{\op} = KQ'/I'$ of type $H$, the corresponding 
cluster-tilted algebra $\G = \End_{\C_H}(T)^{\op}$ is uniquely defined by $Q'$ and $I'$.
\end{cor}

Note that it is well known that a tilted algebra $\L = \End_{H}(T)^{\op}$ can be of finite type, even though
$H$ is not of finite type. 

It would be interesting to know in which generality the main theorem and the above corollary hold.
We do not know any example where it does not hold.

\section{An example}
The path algebra of the quiver 
$$
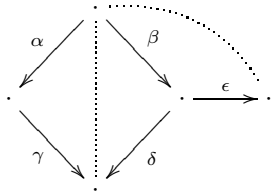
\xy
\xymatrix{
& \cdot \ar[dr]^{\beta} \ar[dl]_{\alpha} \ar@{.}@/^1.1pc/[drr] \ar@{.}[dd] & & \\
\cdot \ar[dr]_{\gamma} & & \cdot \ar[dl]^{\delta} \ar[r]^{\epsilon} & \cdot \\
& \cdot & &
}
\endxy
$$
\\
with relations $\alpha \gamma = \beta \delta$ and $\beta \epsilon = 0$ as indicated,
is a tilted algebra of type $D_5$. The corresponding cluster-tilted algebra is the
path algebra of the quiver
$$
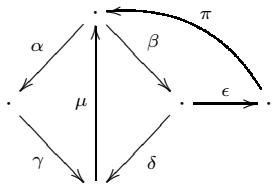
\xy
\xymatrix{
& \cdot \ar[dr]^{\beta} \ar[dl]_{\alpha} & & \\
\cdot \ar[dr]_{\gamma} & & \cdot \ar[dl]^{\delta} \ar[r]^{\epsilon} & \cdot 
\ar@/_1.1pc/[ull]_{\pi} \\
& \cdot \ar[uu]^{\mu}& &
}
\endxy
$$
\\
with relations given by the commutativity relation
$\alpha \gamma = \beta \delta$ and all other compositions of two arrows equal to zero.

\end{document}